\documentclass[12pt,a4paper,reqno]{amsart}
\usepackage{cite}

\usepackage{caption}
\usepackage{subcaption}
\usepackage{amsmath,amssymb,amsfonts,amsthm,epsfig,graphicx,color}
\usepackage{esint}
\usepackage{url}
\DeclareFontFamily{U}{BOONDOX-calo}{\skewchar\font=45 }
\DeclareFontShape{U}{BOONDOX-calo}{m}{n}{
  <-> s*[1.05] BOONDOX-r-calo}{}
\DeclareFontShape{U}{BOONDOX-calo}{b}{n}{
  <-> s*[1.05] BOONDOX-b-calo}{}
\DeclareMathAlphabet{\mathscrboondoxboondox}{U}{BOONDOX-calo}{m}{n}
\SetMathAlphabet{\mathscrboondoxboondox}{bold}{U}{BOONDOX-calo}{b}{n}
\DeclareMathAlphabet{\mathscrboondox}{U}{BOONDOX-calo}{b}{n}

\usepackage{mathrsfs}

\newtheorem{theorem}{Theorem}[section]

\theoremstyle{definition}

\numberwithin{equation}{section}
\numberwithin{figure}{section}

\renewcommand{\le}{\leqslant}

\renewcommand{\ge}{\geqslant}

\def\R{\mathbb R}
\def\e{\varepsilon}
\def\supp{\hbox{supp}}
\def\na{\nabla}
\def\B{\mathcal B}
\newcommand\fb[1]{\partial \{ {#1} > 0 \}}
\title{New trends in free boundary problems}

\author{Serena Dipierro}
\address[Serena Dipierro]{School of Mathematics and Statistics,
University of Melbourne,
813 Swanston Street, Parkville VIC 3010, Australia}
\email{s.dipierro@unimelb.edu.au}

\author{Aram Karakhanyan}
\address[Aram Karakhanyan]{Maxwell Institute for
Mathematical Sciences and School of Mathematics, University of
Edinburgh, James Clerk Maxwell Building, Peter Guthrie Tait Road,
Edinburgh EH9 3FD, United Kingdom}
\email{aram.karakhanyan@ed.ac.uk}

\author{Enrico Valdinoci}
\address[Enrico Valdinoci]{School of Mathematics and Statistics,
University of Melbourne,
813 Swanston Street, Parkville VIC 3010, Australia, 
and Istituto di Matematica Applicata e Tecnologie Informatiche,
Consiglio Nazionale delle Ricerche,
Via Ferrata 1, 27100 Pavia, Italy,
and
Weierstra{\ss} Institut f\"ur Angewandte Analysis und Stochastik,
Mohrenstra{\ss}e 39, 10117 Berlin, Germany,
and Dipartimento di Matematica, Universit\`a degli studi di Milano,
Via Saldini 50, 20133 Milan, Italy}
\email{enrico@mat.uniroma3.it}

\begin{document}

\begin{abstract}
We present a series of recent results
on some new classes of free boundary problems.
Differently from the classical literature, the problems
considered have either a ``nonlocal''
feature (e.g., the interaction or/and the interfacial
energy may
depend on global quantities) or a ``nonlinear''
flavor (namely, the total energy is the nonlinear
superposition of energy components, thus losing
the standard additivity and scale invariances of the problem).

The complete proofs and the full details of the results
presented here are given in~\cite{MR3390089,
MR3427047, Dipierro2016HP, 2015arXiv151203043D, 2016arXiv161100412D, AK}.
\end{abstract}


\maketitle

\begin{center}{\bf Para Don Ireneo, por supuesto.}\end{center}
\bigskip


\section{Introduction}

In this survey, we would like to present some
recent research directions in the study of variational
problems whose minimizers naturally exhibit the formation
of free boundaries. Differently than the cases
considered in most of the existing literature, the problems
that we present here are either {\em nonlinear} (in the sense
that the energy functional is the nonlinear superposition
of classical energy contributions) or {\em nonlocal}
(in the sense that some of the energy contributions
involve objects that depend on the global geometry
of the system).

In these settings, the problems typically show new features
and additional difficulties with respect to the classical cases.
In particular, as we will discuss in further details:
the regularity theory is more complicated,
there is a lack of scale invariance for some problems,
the natural scaling properties of the energy may not
be compatible with the optimal regularity,
the condition at the free boundary may be of nonlocal
or nonlinear type and involve the global behavior of
the solution itself, and 
some problems may exhibit a variational instability
(e.g., minimizers in large domains and in small domains
may dramatically differ the ones from the others).

We will also discuss
how the classical free boundary problems in~\cite{MR618549, MR732100, MR1808651}
are recovered either as limit problems
or after a blow-up, under appropriate structural
conditions on the energy functional.\medskip

We recall the classical free boundary problems
of~\cite{MR618549, MR732100, MR1808651}
in Section~\ref{89102}.

The results concerning nonlocal free boundary problems
will be presented in Section~\ref{SEC2},
while the case of 
nonlinear energy superposition is discussed
in Section~\ref{SEC3}.

\section{Two classical free boundary problems}\label{89102}

A classical problem in fluid dynamics is the description
of a two-dimensional ideal fluid in terms of its stream
function, i.e. of a function whose level sets describe
the trajectories of the fluid.
For this, we consider an incompressible, irrotational
and inviscid fluid which occupies
a given planar region~$\Omega\subset\R^2$.
If~$V:\Omega\to\R^2$ represents the velocity of the particles
of the fluid, the incompressibility condition implies that
the flow of the fluid through any portions of~$\Omega$
is zero (the amount of fluid coming in is exactly the
same as the one going out), that is, 
for any~$\Omega_o\Subset\Omega$, and denoting by~$\nu$
the exterior normal vector,
$$ 0=\int_{\partial\Omega_o} \nabla V\cdot\nu=
\int_{\Omega_o} {\rm div} V.$$
Since this is valid for any subdomain of~$\Omega$, we thus
infer that
\begin{equation}\label{78:SIDIC}
{\rm div} V=0\; {\mbox{ in }}\Omega.
\end{equation}
Now, we use that the fluid is irrotational to write equation~\eqref{78:SIDIC}
as a second order PDE. To this aim, let us analyze what a ``vortex'' is.
Roughly speaking, a vortex is given by a close trajectory,
say~$\gamma: S^1\to\Omega$,
along which the fluid particles move. In this way, the velocity
field~$V$ is always parallel to the tangent direction~$\gamma'$
and therefore
\begin{equation}\label{78:SIDIC2} 0\ne \int_{S^1} V(\gamma(t))\cdot\gamma'(t)\,dt
=\oint_\gamma V,\end{equation}
where the standard notation for the
circulation line integral is used.
That is, if we denote by~$S$ the region inside~$\Omega$
enclosed by the curve~$\gamma$
(hence, $\gamma=\partial S$), 
we infer by~\eqref{78:SIDIC2} and Stokes' Theorem that
\begin{equation}\label{78:SIDIC212} 0\ne 
\int_S {\rm curl} V\cdot e_3,\end{equation}
where, as usual, we write~$\left\{e_1,e_2,e_3\right\}$
to denote the standard basis of~$\R^3$,
we identify the vector~$V=(V_1,V_2)$
with its three-dimensional image~$V=(V_1,V_2,0)$, and
$$ {\rm curl} V(x):= \nabla \times V (x)=
\det \left( 
\begin{matrix}
e_1 & e_2 & e_3\\
\partial_{x_1}&\partial_{x_2}&\partial_{x_3}\\
V_1(x_1,x_2)&V_2(x_1,x_2)&0
\end{matrix}
\right)=\left( \partial_{x_1}V_2(x)-\partial_{x_2}V_1(x)\right)\,e_3.
$$
In this setting, the fact that the fluid is irrotational is translated
in mathematical language into the fact that the opposite of~\eqref{78:SIDIC212}
holds true, namely
\begin{equation*} 0=
\int_S {\rm curl} V\cdot e_3,\end{equation*}
for any~$S\subset\Omega$ (say, with smooth boundary).
Since this is valid for any arbitrary region~$S$, we thus
can translate the irrotational property of the fluid into 
the condition~${\rm curl} V=0$ in~$\Omega$, that is
\begin{equation}\label{PIRRO}
\partial_{x_1}V_2-\partial_{x_2}V_1=0\;{\mbox{ in }}\Omega.
\end{equation}
Now, we consider the $1$-form 
\begin{equation}\label{OM78:A}\omega:=
V_2\,dx_1-V_1\,dx_2,\end{equation} 
and we have that
$$ d\omega= -\partial_{x_2}V_2\,dx_1\wedge dx_2-
\partial_{x_1}V_1\,dx_1\wedge dx_2
=-{\rm div} V\,dx_1\wedge dx_2
=0,$$ thanks to~\eqref{78:SIDIC}.
Namely $\omega$ is closed, and thus exact (by 
Poincar\'e Lemma, at least if~$\Omega$ is star-shaped).
This says that there exists a function~$u$ such that
$$\omega=du=\partial_{x_1}u\,dx_1+
\partial_{x_2}u\,dx_2.$$
By comparing this and~\eqref{OM78:A}, we conclude that
\begin{equation} \label{67q	8wsIPAa}\partial_{x_1}u=V_2\;
{\mbox{ and }}\;\partial_{x_2}u=-V_1.\end{equation}
We observe that $u$ is a stream function for the fluid, namely
the fluid particles move along the level sets of~$u$:
indeed, if~$x(t)$ is the position of the fluid particle at time~$t$,
we have that~$\dot x(t)=V(x(t))$ is the velocity of the fluid, and
\begin{eqnarray*} &&
\frac{d}{dt} u(x(t))=\partial_{x_1} u(x(t))\dot x_1(t)+
\partial_{x_2} u(x(t))\dot x_2(t)\\&&\quad
=
\partial_{x_1} u(x(t))V_1(x(t))+
\partial_{x_2} u(x(t))V_2(x(t))
=V_2(x(t))V_1(x(t))-V_1(x(t))V_2(x(t))
=0,\end{eqnarray*}
in view of~\eqref{67q	8wsIPAa}.

The stream function~$u$ also satisfies a natural overdetermined
problem. First of all, since~$\partial\Omega$ represents the boundary
of the fluid, and the fluid motion occurs on the level sets of~$u$,
up to constants we may assume that~$u=0$ along~$\partial\Omega$.
In addition, along~$\partial\Omega$ Bernoulli's Law prescribes
that the velocity is balanced by the pressure (which we take
here to be $p=p(x)$). That is,
up to dimensional constants, we can write that, along~$\partial\Omega$,
$$ p = |V|^2 =|\nabla u|^2,$$
where we used again~\eqref{67q	8wsIPAa} in the last identity.

Also, \eqref{PIRRO} and~\eqref{67q	8wsIPAa} give that, in~$\Omega$,
$$ \Delta u=\partial_{x_1}V_2-\partial_{x_2}V_1=0,$$
that is, summarizing,
\begin{equation} \label{78:APAJ}\left\{
\begin{matrix}
\Delta u=0 & {\mbox{ in }}\Omega,\\
u=0& {\mbox{ on }}\partial\Omega,\\
|\nabla u|^2=p& {\mbox{ on }}\partial\Omega.
\end{matrix}
\right. \end{equation}
Notice that these types of overdetermined problems
are, in general, 
not solvable: namely, only ``very special'' domains
allow a solution of such overdetermined problem to exist
(see e.g.~\cite{MR0333220}). In this spirit, determining such
domain~$\Omega$ is part of the problem itself, and
the boundary of~$\Omega$ is, in this sense, a ``free boundary''
to be determined together with the solution~$u$.

These kinds of free boundary problems
have a natural formulation, which was widely studied
in~\cite{MR618549, MR732100}.
The idea is to consider an energy functional
which is the superposition of a Dirichlet part
and a volume term. By an appropriate domain variation,
one sees that minimizers (or, more generally, critical
points) of this functional correspond (at least in a weak sense)
to solutions of~\eqref{78:APAJ} (compare, for instance,
the system in~\eqref{78:APAJ}
here with Lemma~2.4 and Theorem~2.5 in~\cite{MR618549}).
Needless to say, in this framework, the analysis of the minimizers
of this energy functional and of their level sets
becomes a crucial topic of research. 
\medskip

In~\cite{MR1808651} a different energy functional
is taken into account, in which the volume term
is substituted by a perimeter term. This modification provides
a natural change in the free boundary condition (in this
setting, the pressure of the Bernoulli's Law is replaced
by the curvature of the level set, see formula~(6.1)
in~\cite{MR1808651}).\medskip

In the following sections we will discuss what happens when:
\begin{itemize}
\item we interpolate the volume term of
the energy functional of~\cite{MR618549, MR732100}
and the perimeter term of
the energy functional of~\cite{MR1808651}
with a fractional perimeter term, which recovers the volume
and the classical perimeter in the limit;
\item we consider a nonlinear energy superposition,
in which the total energy depends on the volume, or on the
(possibly fractional) perimeter, in a nonlinear fashion.
\end{itemize}

\section{Nonlocal free boundary problems}\label{SEC2}

A classical motivation for free boundary problems
comes from the superposition of
a ``Dirichlet-type energy''~${\mathscrboondox{D}}$ 
and an ``interfacial energy''~${\mathscrboondox{I}}$.
Roughly speaking, one may consider the minimization problem
of an energy functional
\begin{equation} \label{ENERGIA LIN}
{\mathscrboondox{E}}:={\mathscrboondox{D}}+{\mathscrboondox{I}},\end{equation}
which takes into account the following two 
tendencies of the energy contributions, namely:
\begin{itemize}
\item the term ${\mathscrboondox{D}}$
tries to reduce the oscillations of the minimizers,
\item
while the term~${\mathscrboondox{I}}$ penalizes the formation
of interfaces.\end{itemize}

Two classical approaches
appear in the literature to measure these interfaces,
taking into account the ``volume'' of the phases
or the ``perimeter'' of the phase separations.
The first approach, based on a ``bulk'' energy contribution,
was widely studied in~\cite{MR618549}.
In this setting, the energy superposition
in~\eqref{ENERGIA LIN} (with respect to
a reference domain~$\Omega\subset\R^n$) takes the form
\begin{equation}\label{AC-for}
\begin{split}
& {\mathscrboondox{D}}={\mathscrboondox{D}}(u):=\int_\Omega |\nabla u(x)|^2\,dx
\\{\mbox{and}}\quad&{\mathscrboondox{I}}={\mathscrboondox{I}}(u):=
\int_\Omega \chi_{ \{ u>0\} }(x)\,dx
= {\mathscr{L}}^n \big( \Omega\cap \{ u>0\} \big),\end{split}\end{equation}
where~${\mathscr{L}}^n$ denotes, as customary,
the $n$-dimensional Lebesgue measure.
The case of two phase contributions
(namely, the one which takes into
account
the bulk energy of both $\{u>0\}$ and~$\{u<0\}$)
was also considered in~\cite{MR732100}.
\medskip

The second approach, based on a ``surface tension''
energy contribution,
was introduced in~\cite{MR1808651}.
In this setting, the energy superposition
in~\eqref{ENERGIA LIN} takes the form
\begin{equation}\label{ACKS-for}
\begin{split}
& {\mathscrboondox{D}}={\mathscrboondox{D}}(u):=\int_\Omega |\nabla u(x)|^2\,dx
\\{\mbox{and}}\quad&{\mathscrboondox{I}}={\mathscrboondox{I}}(u):=
{\rm Per}\big(\{ u>0\} ,\Omega\big),\end{split}\end{equation}
where the notation
$${\rm Per}(E,\Omega):=\int_\Omega |D\chi_E(x)|\,dx
= [ \chi_E]_{BV(\Omega)}$$
represents the perimeter of the set~$E$ in~$\Omega$;
hence, if~$E$ has smooth boundary, then
${\rm Per}(E,\Omega)={\mathscr{H}}^{n-1}\big((\partial E)\cap\Omega\big)$,
being~${\mathscr{H}}^{n-1}$ the $(n-1)$-dimensional
Hausdorff measure.\medskip

As pointed out in~\cite{MR3390089},
the two free boundary problems
in~\eqref{AC-for} and~\eqref{ACKS-for}
can be settled into a unified framework,
and in fact they may be seen as ``extremal'' problems
of a family of energy functionals
indexed by a continuous parameter~$\sigma\in(0,1)$.

To this aim, given two measurable
sets~$E$, $F\subset\R^n$, with~${\mathscr{L}}^n(E\cap F)=0$,
one considers the $\sigma$-interaction of~$E$ and~$F$,
as given by the double integral
$$ {\mathscr{S}}_\sigma(E,F):=
\sigma\,(1-\sigma)\,\iint_{E\times F}\frac{dx\,dy}{|x-y|^{n+\sigma}}.$$
In~\cite{MR2675483}, the notion
of $\sigma$-minimal surfaces has been introduced
by considering minimizers of the $\sigma$-perimeter
induced  by such interaction. Namely, one defines
the $\sigma$-perimeter of~$E$ in~$\Omega$
as the contribution relative to~$\Omega$
of the $\sigma$-interaction of~$E$ and its complement (which
we denote by~$E^c:=\R^n\setminus E$), that is
\begin{equation}\label{SIGP} {\rm Per}_\sigma(E,\Omega)
:= {\mathscr{S}}_\sigma(E,E^c\cap \Omega)+
{\mathscr{S}}_\sigma(E\cap\Omega,E^c\cap\Omega^c).\end{equation}
After~\cite{MR2675483},
an intense activity has been performed to investigate
the regularity and the geometric properties
of $\sigma$-minimal surfaces:
see in particular~\cite{MR3090533, MR3035063, MR3107529,
FIGALLI, MR3331523}
for interior regularity results, and~\cite{MR3532394, MR3516886, Dipierro2016JFA, 2016arXiv161208295B}
for the rather special behavior of
$\sigma$-minimal surfaces near the boundary
of the domain. See also~\cite{2016arXiv160706872D}
for a recent survey on $\sigma$-minimal surfaces.
\medskip

The analysis of the asymptotics of the $\sigma$-perimeter
as~$\sigma\nearrow1$ has been studied
under several perspectives in~\cite{MR1945278, MR1942130, MR2033060,
MR2782803, MR2765717, MR3107529, MR3556344}. Roughly speaking,
up to normalization constants, we may say that~${\rm Per}_\sigma$
approaches the classical perimeter as~$\sigma\nearrow1$.
On the other hand, as~$\sigma\searrow0$,
we have that~${\rm Per}_\sigma$
approaches the Lebesgue measure (again, up
to normalization constants,
see~\cite{MR1940355, MR3007726, MR3506705}
for precise statements and examples).

In virtue of these considerations,
we have that the free boundary problem
introduced in~\cite{MR3390089},
which takes into account
the energy superposition
in~\eqref{ENERGIA LIN} of the form
\begin{equation}\label{CSV-for}
\begin{split}
& {\mathscrboondox{D}}={\mathscrboondox{D}}(u):=\int_\Omega |\nabla u(x)|^2\,dx
\\{\mbox{and}}\quad&{\mathscrboondox{I}}={\mathscrboondox{I}}(u):=
{\rm Per}_\sigma\big(\{ u>0\} ,\Omega\big),\end{split}\end{equation}
may be seen as an interpolation
of the problems stated in~\eqref{AC-for} and~\eqref{ACKS-for}
(that is, at least at a formal level, 
the energy functional in~\eqref{CSV-for} reduces to that
in~\eqref{AC-for} as~$\sigma\searrow0$
and to that in~\eqref{ACKS-for}
as~$\sigma\nearrow1$). \medskip

A nonlocal variation of the classical Dirichlet energy
has been also considered in~\cite{MR2677613, MR2926238}.
In this setting, the classical $H^1$-seminorm
in~$\Omega$ of a function~$u$ is replaced
by a Gagliardo $H^s$-seminorm of the form
\begin{equation} \label{GAGLIARDO} s\,(1-s)\,
\iint_{Q_\Omega} \frac{|u(x)-u(y)|^2}{|x-y|^{n+2s}}\,dx\,dy,\end{equation}
where
\begin{equation}\label{NOT Q}
Q_\Omega:=(\Omega\times\Omega)\cup
(\Omega\times\Omega^c)\cup (\Omega^c\times\Omega)\end{equation}
and~$s\in(0,1)$. More precisely,
in~\cite{MR2677613, MR2926238} a superposition
of the Gagliardo seminorm
and the Lebesgue measure of the positivity
set is taken into account.

It is worth to point out that the domain~$Q_\Omega$
in~\eqref{GAGLIARDO} comprises all the interactions
of points~$(x,y)\in\R^{2n}$ which involve the domain~$\Omega$,
since~$ Q_\Omega=(\R^{n}\times\R^n)\setminus 
(\Omega^c\times\Omega^c)$.

In this sense, the integration over~$Q_\Omega\subset\R^{2n}$
is the natural counterpart of the classical integration over~$\Omega$
of the standard Dirichlet energy, since~$\Omega=\R^n\setminus\Omega^c$.

Also, the double integral in~\eqref{GAGLIARDO}
recovers the classical Dirichlet energy, see e.g.~\cite{MR1945278, MR2944369}.
\medskip

In this spirit, in~\cite{MR3427047, Dipierro2016HP}
a fully nonlocal counterpart of the free boundary
problems in~\eqref{AC-for} and~\eqref{ACKS-for}
has been introduced, by studying
energy superpositions of Gagliardo norms
and fractional perimeters (see also~\cite{LUCA}
for the 
superpositions of Gagliardo norms
and classical perimeters).
More precisely, 
the energy superposition
in~\eqref{ENERGIA LIN}
considered in~\cite{MR3427047, Dipierro2016HP} takes the form
\begin{equation*}
\begin{split}
& {\mathscrboondox{D}}={\mathscrboondox{D}}(u):=
s\,(1-s)\,
\iint_{Q_\Omega} \frac{|u(x)-u(y)|^2}{|x-y|^{n+2s}}\,dx\,dy
\\{\mbox{and}}\quad&{\mathscrboondox{I}}={\mathscrboondox{I}}(u):=
{\rm Per}_\sigma\big(\{ u>0\} ,\Omega\big),\end{split}\end{equation*}
where~$s$, $\sigma\in(0,1)$.\medskip

We summarize here a series of results recently obtained
in~\cite{MR3390089, MR3427047, Dipierro2016HP}
for these nonlocal free boundary problems
(some of these results also rely on a notion of fractional harmonic
replacement analyzed in~\cite{MR3320130}).
First of all, we have that
minimizers\footnote{Here,
for simplicity,
we omit the fact that, in this setting, the minimization
is performed not only on a function, but
on a couple given by the function and its
positivity set. See Section~2 in~\cite{MR3390089}
for a rigorous discussion on this important
detail.}
of free boundary problems
with fractional perimeter interfaces
are continuous, possess suitable density estimates
and have smooth free boundaries
up to sets
of codimension~3:

\begin{theorem}\label{CC67}
[Theorems 1.1 and 1.2 in \cite{MR3390089}]
Let~$u_\star$ be a minimizer of
\begin{equation}\label{EN+CSV} 
{\mathscrboondox{E}}(u):=\int_{B_1} |\nabla u(x)|^2\,dx+
{\rm Per}_\sigma\big( \{u>0\},B_1\big),\end{equation}
with~$\sigma\in(0,1)$ and~$0\in\partial\{u_\star>0\}$.

Then~$ u_\star$ is locally~$C^{1- \frac\sigma2}$
and, for any~$r\in\left(0,\frac12\right)$,
$$￼ \min\left\{ {\mathscr{L}}^n(B_r\cap \{u_\star\ge0\}),\,
{\mathscr{L}}^n(B_r\cap \{u_\star\le0\}) \right\}\ge cr^n,$$
for some~$c>0$.

Moreover, the free boundary is a $C^{\infty}$-hypersurface
possibly outside a small singular
set of Haussdorff dimension~$n-3$.
\end{theorem}

We remark that the H\"older exponent~$1- \frac\sigma2$
is consistent with the natural scaling of the problem (namely~$u_r
(x) := r^{\frac\sigma2-1} u_\star (rx)$ is still a minimizer).
Such type of regularity approaches the optimal exponent
in~\cite{MR618549, MR732100}
as~$\sigma\searrow0$. Nevertheless, as~$\sigma\nearrow1$,
minimizers in~\cite{MR1808651} are known to be Lipschitz
continuous, therefore we think that it is a very interesting
open problem to investigate the optimal
regularity of~$u_\star$ in Theorem~\ref{CC67}
(we stress that this optimal regularity may well approach
the Lipschitz regularity and so ``beat the natural scaling
of the problem'').

Also, we think it is very interesting to obtain
optimal bounds on the dimension of the singular set
in Theorem~\ref{CC67}.\medskip

It is also worth to observe that the minimizers
in Theorem~\ref{CC67} satisfy a nonlocal free boundary
condition. Namely,
the normal jump~${J}_\star:=
|\nabla u_\star^+|^2-|\nabla u_\star^-|^2$
along the smooth points
of the free boundary coincides (up to normalizing constants)
with the nonlocal mean curvature of the free boundary,
which is defined by
\begin{equation}\label{KADEF} 
{\mathscr{K}}^\sigma(x):=\int_{\R^n} \frac{\chi_{\{ u\le0\}}(y)-\chi_{\{ u>0\}}(y)}{|x-y|^{n+\sigma}}\,dy,\end{equation}
for $x\in\partial\{ u>0\}$.

This free boundary condition has been presented in formula~(1.6)
of~\cite{MR3390089}. Since~${\mathscr{K}}^\sigma$ approaches
the classical mean curvature as~$\sigma\nearrow1$
and a constant as~$\sigma\searrow0$ (see e.g.~\cite{MR3230079}
and
Appendix~B 
of~\cite{2016arXiv160706872D}), we remark
that this nonlocal free boundary condition
recovers the classical ones in~\cite{MR1808651}
and in~\cite{MR618549, MR732100} as~$\sigma\nearrow1$,
and as~$\sigma\searrow0$, respectively.
\medskip

In~\cite{MR3427047, Dipierro2016HP},
we consider the fully nonlocal case in which both the
energy components
become of nonlocal type,
namely we replace~\eqref{EN+CSV} with the energy
functional
\begin{equation}\label{EN+FUL} 
{\mathscrboondox{E}}(u):=
s\,(1-s)\,
\iint_{Q_{B_1}} \frac{|u(x)-u(y)|^2}{|x-y|^{n+2s}}\,dx\,dy
+
{\rm Per}_\sigma\big( \{u>0\},B_1\big),\end{equation}
with~$s$, $\sigma\in(0,1)$, where the notation in~\eqref{NOT Q}
has been also used.

In this setting, we have:

\begin{theorem}\label{09ijzaxis}[Theorem 1.1 in~\cite{Dipierro2016HP}]
Let~$u_\star$ be a minimizer of~\eqref{EN+FUL} 
with~$0\in\partial\{u_\star>0\}$.

Assume that~$u_\star\ge0$ in~$B_1^c$ and that
$$ \int_{\R^n} \frac{|u_\star(x)|}{1+|x|^{n+2s}}\,dx<+\infty.$$
Then,~$ u_\star$ is locally~$C^{s- \frac\sigma2}$
and, for any~$r\in\left(0,\frac12\right)$,
$$￼ \min\left\{ {\mathscr{L}}^n(B_r\cap \{u_\star\ge0\}),\,
{\mathscr{L}}^n(B_r\cap \{u_\star=0\}) \right\}\ge cr^n,$$
for some~$c>0$.
\end{theorem}

We observe that the H\"older exponent in Theorem~\ref{09ijzaxis}
recovers that of Theorem~\ref{CC67} as~$s\nearrow1$.
Once again, we think that it would be very interesting
to investigate the optimal regularity of the minimizers
in Theorem~\ref{09ijzaxis}. Also,
Theorem~\ref{09ijzaxis} has been established in the
``one-phase'' case, i.e. under the assumption that the
minimizer has a sign. It would be very interesting
to establish similar results
in the ``two-phase'' case in which minimizers can change
sign. It is worth to remark that the case in which
minimizers change sign is conceptually
harder in the nonlocal setting
than in the local one, since the two phases interact between
each other, thus producing additional energy contributions
which need to be carefully taken into account.

\section{Nonlinear free boundary problems}\label{SEC3}

In~\cite{2015arXiv151203043D, 2016arXiv161100412D}
a new class of free boundary problems
has been considered, by taking into account
``nonlinear energy superpositions''.
Namely, differently than in~\eqref{ENERGIA LIN},
the total energy functional considered in~\cite{2015arXiv151203043D, 2016arXiv161100412D}
is of the form
\begin{equation} \label{ENERGIA NON-LIN}
{\mathscrboondox{E}}:={\mathscrboondox{D}}+\Phi_0({\mathscrboondox{I}}),\end{equation}
for a suitable function~$\Phi_0$.
When~$\Phi_0$ is linear, the energy functional in~\eqref{ENERGIA NON-LIN}
boils down to its ``linear counterpart''
given in~\eqref{ENERGIA LIN}, but for a
nonlinear function~$\Phi_0$
the minimizers of the energy functional in~\eqref{ENERGIA NON-LIN}
may exhibit\footnote{Let us give a brief motivation
for the nonlinear interface case. The classical energy functionals
in~\cite{MR618549, MR732100, MR1808651}
may be also considered in view of models arising in population 
dynamics.
Namely, one can consider the regions~$\{u>0\}$ and~$\{u\le0\}$
as areas occupied by two different populations, that have
reciprocal ``hostile'' feelings. Then, the diffusive behavior
of the populations (which is
encoded by the Dirichlet term of the energy)
is influenced by the fact that the two populations will have the
tendency to minimize the contact between themselves,
and so to reduce an interfacial energy as much as possible.

In this setting, it is natural to consider the case
in which the reaction of the populations to the mutual contact
occurs in a nonlinear way. For instance, the case in which
additional irrationally motivated hostile feelings arise
from further contacts between the populations is naturally
modeled by a superlinear interfacial energy, while the
case in which the interactions between the populations
favor the possibility of compromises and cultural exchanges
is naturally
modeled by a sublinear interfacial energy.}
interesting differences with respect to the classical
case.\medskip

A detailed analysis of free boundary problems as in~\eqref{ENERGIA NON-LIN}
is given in~\cite{2015arXiv151203043D, 2016arXiv161100412D}.
Here, we summarize some of the results obtained
(we give here simpler statements, referring
to~\cite{2015arXiv151203043D, 2016arXiv161100412D}
for more general results).
We take~$\Phi_0:[0,+\infty)\to[0,+\infty)$ to
be monotone, nondecreasing,
lower semicontinuous and coercive -- in the sense that
$$ \lim_{t\nearrow+\infty} \Phi_0(t)=+\infty.$$
We will also use the notation of writing~${\rm Per}_\sigma$
for every~$\sigma\in[0,1]$, with the convention that
\begin{itemize}
\item when~$\sigma\in(0,1)$, ${\rm Per}_\sigma$
is the nonlocal perimeter defined in~\eqref{SIGP};
\item when~$\sigma=1$, ${\rm Per}_\sigma$ is the classical
perimeter;
\item when~$\sigma=0$, ${\rm Per}_\sigma(E;\Omega)=
{\mathscr{L}}^n(E\cap\Omega)$.
\end{itemize}
Then, in the spirit of~\eqref{ENERGIA NON-LIN}, we consider energy functionals of the form
\begin{equation}\label{NLI}
{\mathscrboondox{E}}(u):=\int_{\Omega} |\nabla u(x)|^2\,dx+\Phi_0\left(
{\rm Per}_\sigma\big( \{u>0\},\Omega\big)\right).\end{equation}
Notice that, for~$\sigma\in(0,1)$ and~$\Phi_0(t)=t$,
the energy in~\eqref{NLI} reduces to that in~\eqref{EN+CSV}.
Similarly, for~$\sigma=0$ and~$\sigma=1$,
the energy in~\eqref{NLI}
boils down 
to those in~\cite{MR618549} and~\cite{MR1808651}, respectively.

When~$\sigma=0$,
a particularly interesting case of nonlinearity is given by~$\Phi_0(t)=
t^{\frac{n-1}{n}}$. Indeed, in this case, the interfacial energy
depends on the $n$-dimensional Lebesgue measure,
but it scales like an $(n-1)$-dimensional surface measure
(also, by Isoperimetric Inequality, the energy levels
of the functional in~\cite{MR1808651}
are above those in~\eqref{NLI}).\medskip

We point out that the free boundary problems in~\eqref{NLI}
develop a sort of natural instability, in the sense that
minimizers in a large ball, when restricted to smaller balls,
may lose their minimizing properties. In fact, minimizers
in large and small balls may be rather different from each
other:

\begin{theorem}\label{890AJ}[Theorem~1.1 
in~\cite{2015arXiv151203043D}]
There exist
a nonlinearity~$\Phi_0$
and radii~$R_0>r_0>0$ such that
a minimizer
for~\eqref{NLI} in~$\Omega:=B_{R_0}$
is not a minimizer
for~\eqref{NLI} in~$\Omega:=B_{r_0}$.
\end{theorem}

The counterexample in Theorem~\ref{890AJ}, which clearly shows the
lack\footnote{Just to recall the importance of scaling invariances in the classical free
boundary problems, let us quote page 114 of \cite{MR732100}:
``for (small) balls $B_r$ [...] let us assume (see 3.1) $B_r=B_1(0)$''.}
of scaling invariance of the problem,
is constructed by taking advantage of the different
rates of scaling produced by a suitable nonlinear function~$\Phi_0$,
chosen to be constant on an interval.
Namely, the ``saddle function'' in the plane
$u_0(x_1,x_2)=x_1x_2$ is harmonic and therefore
minimizes the Dirichlet energy of~\eqref{NLI}.
In large balls, the interface of~$u_0$ (as well as the ones
of its competitors) produces a contribution that
lies in the constant part of~$\Phi_0$, thus reducing
the minimization problem of~\eqref{NLI} to the one
coming from the Dirichlet contribution, and so favoring~$u_0$
itself. Viceversa, in small balls, the interface
of~$u_0$ produces more (possibly fractional) perimeter than
the one of the competitors whose positivity sets do not
come to the origin, and this fact implies that~$u_0$
is not a minimizer in small balls.


\begin{figure}
        \centering
                \begin{subfigure}[b]{0.40\textwidth}
                                \includegraphics[width=\textwidth]{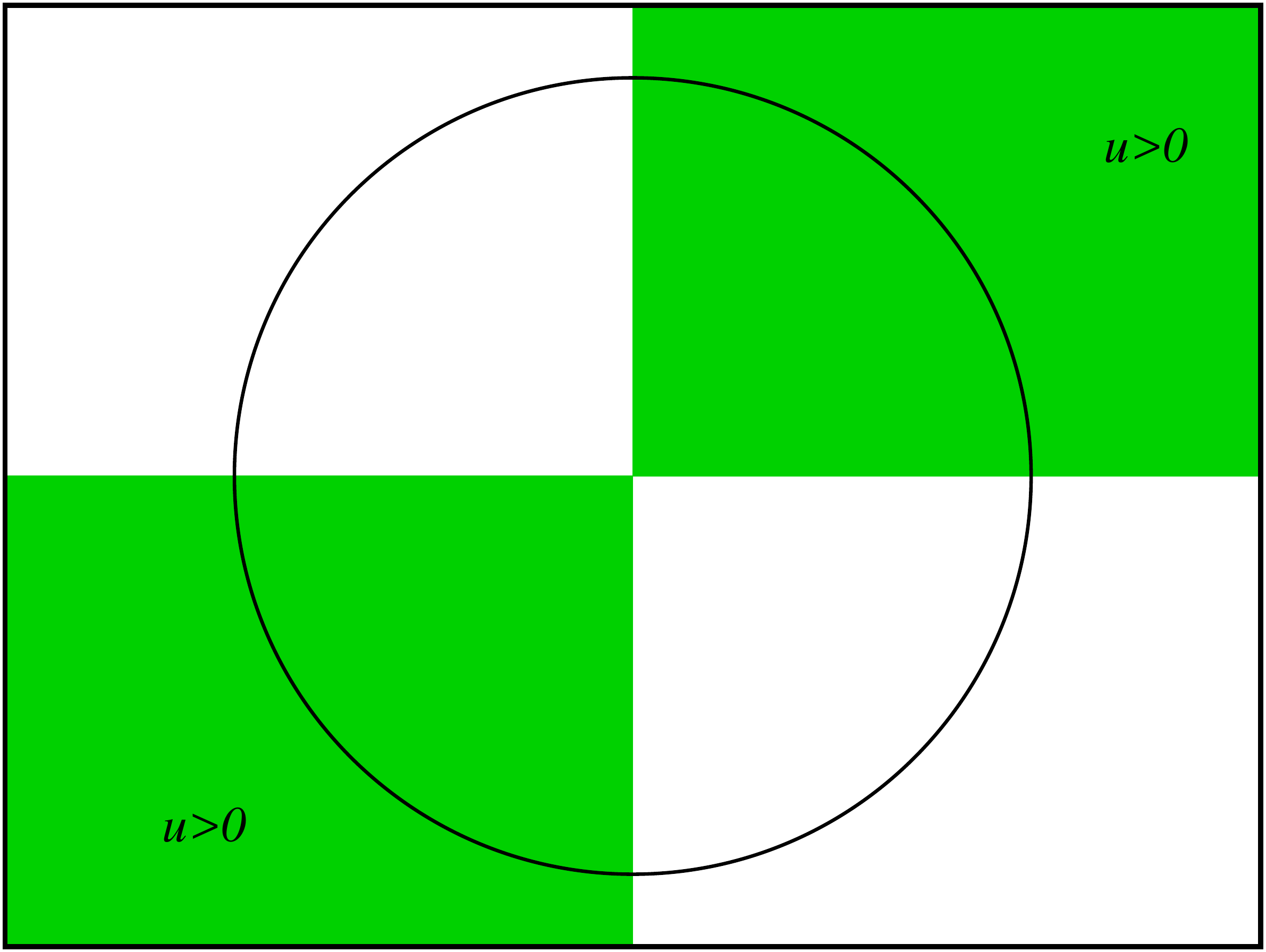}
                \caption{The minimizer in Theorem \ref{890AJ}
for a large ball $B_{R_0}$.}
                       \end{subfigure}
                       \hspace{1cm}
                \begin{subfigure}[b]{0.40\textwidth}
                \includegraphics[width=\textwidth]{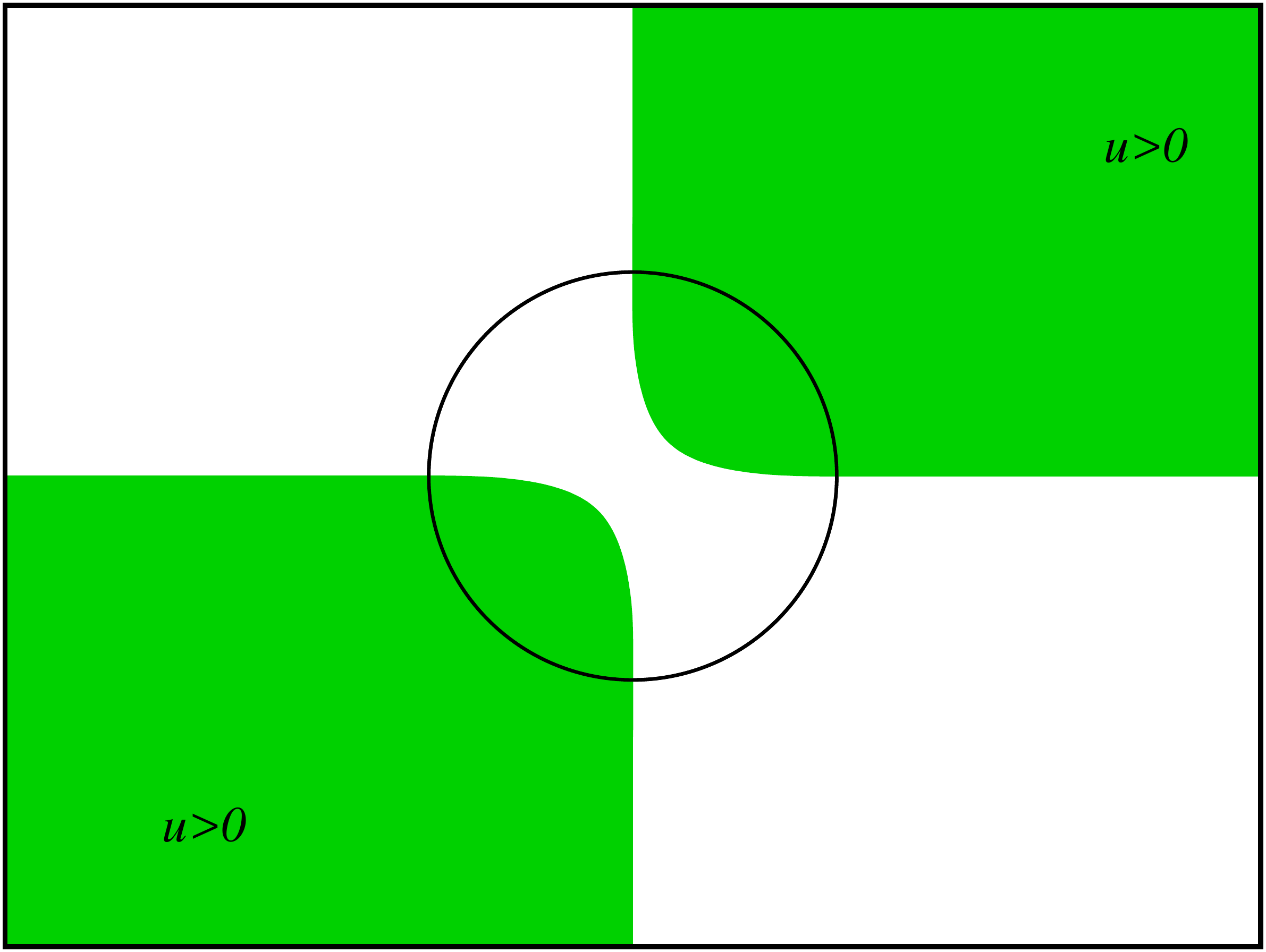}
                \caption{A minimizer in Theorem \ref{890AJ}
for a small ball $B_{r_0}$.}
        \end{subfigure}
              \caption{}
              \label{NESS}
                \end{figure}

This argument, which is depicted in Figures~\ref{NESS} (A) and (B),
is rigorously explained in Section~3
of~\cite{2015arXiv151203043D}.
\medskip

In spite of this instability and of the lack of self-similar properties
of the energy functional, some regularity results
for minimizers of~\eqref{NLI} still hold true, under appropriate
assumptions on the nonlinearity~$\Phi_0$
(notice that, for~$\Phi_0$ with a constant part,
Theorem~\ref{890AJ} would produce the minimizer~$u_0$
whose free boundary is a singular
cone, see Figure~\ref{NESS} (A), hence any regularity
result on the free boundary has to rule 
out this possibility by a suitable assumption on~$\Phi_0$).
In this sense, we have the following results:

\begin{theorem}\label{18ueyy23r}
[Corollary 1.4 and Theorems 1.5 and 1.6
in~\cite{2015arXiv151203043D}]
Let~$\sigma\in(0,1]$, $\Phi_0$ be Lipschitz continuous
and strictly increasing.
Let~$u_\star$ be a minimizer of~\eqref{NLI}
in~$\Omega:=B_R$, with~$0\le u_\star\le M$
on~$\partial B_R$, for some~$M>0$.

Then, $u_\star\in C^{1-\frac\sigma2}(B_{R/4})$.
Also, for any~$r\in\left(0,\frac{R}{4}\right)$,
$$￼ \min\left\{ {\mathscr{L}}^n(B_r\cap \{u_\star\ge0\}),\,
{\mathscr{L}}^n(B_r\cap \{u_\star=0\}) \right\}\ge cr^n,$$
for some~$c>0$.
\end{theorem}

For~$\sigma=0$, a result similar to that in Theorem~\ref{18ueyy23r}
holds true, in
the sense that~$u_\star$ is Lipschitz,
see Theorems~1.3, 8.1 and~9.2 in~\cite{2016arXiv161100412D}.
Moreover, in this case one obtains additional results, such as the
nondegeneracy of the minimizers,
the partial regularity of the free boundary and the
full regularity in the plane:

\begin{theorem}\label{18ueyy23r:BIS}
[Theorems 1.4, 1.6 and 1.7
in~\cite{2016arXiv161100412D}]
Let~$\sigma=0$, $\Phi_0$ be Lipschitz continuous
and strictly increasing.
Let~$u_\star$ be a minimizer of~\eqref{NLI}
in~$\Omega$, with~$0\in\partial\{u_\star>0\}$
(in the measure theoretic sense). 

Then, for any~$D\Subset\Omega$,
there exists~$c>0$ such that for any~$r>0$ for which~$B_r\Subset D$,
it holds that
$$ \int_{B_r\cap\{u_\star>0\} } u_\star^2(x)\,dx\ge c r^{n+2}.$$
Also, $\nabla u_\star$ is locally BMO, in the sense that
$$ \sup_{B_r\Subset D}\fint_{B_r} \big| \nabla u_\star(x)-
\langle \nabla u_\star \rangle_r \big|\,dx
\le C,$$
for some~$C>0$, where
$$ \langle \nabla u_\star \rangle_r:=\fint_{B_r} \nabla u_\star (x)\,dx.$$
In addition~${\mathcal{H}}^{n-1}\big(B_r\cap
(\partial\{u_\star>0\})\big) <+\infty$.

Finally, if~$n=2$, then~$B_r\cap
(\partial\{u_\star>0\})$ is made of continuously differentiable curves.
\end{theorem}

The BMO-type regularity and the partial regularity
of the free boundary in Theorem~\ref{18ueyy23r:BIS}
rely in turn on some techniques developed in~\cite{2015arXiv150807447D}.

It is also interesting to remark that the case~$\sigma=0$
recovers the classical problems in~\cite{MR618549}
after a blow-up:

\begin{theorem}\label{18ueyy23r:TRIS}
[Theorem 1.5 and Proposition~10.1
in~\cite{2016arXiv161100412D}]
Let~$\sigma=0$, $\Phi_0$ be Lipschitz continuous
and strictly increasing.
Let~$u_\star$ be a minimizer of~\eqref{NLI}
in~$\Omega$, with~$0\in\Omega$. 

For any~$r>0$, let~$u_r(x):=
\frac{ u_\star (rx)}{r}$.
Then, there exists the blow-up limit
$$ u_0(x):=\lim_{r\searrow0} u_r(x).$$
Also, $u_0$ is continuous and with linear growth,
and it is a minimizer of the functional
\begin{equation}\label{E0AL} {\mathscrboondox{E}}_0(u):=
\int_{B_\rho} |\nabla u(x)|^2\,dx+
\lambda_0\,{\mathcal{L}}^n\big(B_\rho\cap\{u>0\}\big),\end{equation}
where
\begin{equation}\label{E0AL2} \lambda_0:= \Phi_0'\left(
{\mathcal{L}}^n\big(\Omega\cap\{u_\star>0\}\big)
\right).\end{equation}
\end{theorem}

We stress that the energy functional in~\eqref{E0AL}
coincides with that analyzed in the classical paper~\cite{MR618549}.
Nevertheless, the ``scaling constant''~$\lambda_0$
in~\eqref{E0AL} depends on the original minimizer~$u_\star$,
as prescribed by~\eqref{E0AL2} (only in the case of a
linear~$\Phi_0$, we have that~$\lambda_0$
is a structural constant independent of~$u_\star$).\medskip

The fact that geometric and physical quantities
arising in this type of problem are not universal constants
but depend on the minimizer itself
is, in our opinion, an intriguing feature of this type of
problems. In this sense, we recall that in~\cite{MR618549}
the free boundary condition coincides with the classical
Bernoulli's law, namely
the normal jump~${J}_\star:=
|\nabla u_\star^+|^2-|\nabla u_\star^-|^2$
along the smooth points
of the free boundary is constant
(in~\cite{MR1808651}
it coincides with the mean curvature
of the free boundary).
Differently from the classical cases,
in our nonlinear setting, the free boundary
condition depends on the minimizer itself.
Indeed, 
in our case the normal jump~${J}_\star$ coincides with
\begin{equation}\label{78:0ass12d}
{\mathscr{K}}^\sigma\;\Phi_0'\left( 
{\rm Per}_\sigma(\{u_\star>0\},\Omega)\right),
\end{equation}
where~${\mathscr{K}}^\sigma$ is the 
nonlocal mean curvature of the free boundary, as defined
in~\eqref{KADEF} 
(see formula~(1.12) in~\cite{2015arXiv151203043D}
and formulas~(1.13) and~(1.14) in~\cite{2016arXiv161100412D}).

We point out that~\eqref{78:0ass12d}
recovers the classical cases in~\cite{MR618549, MR1808651}
when~$\sigma\in\{0,1\}$ and~$\Phi_0$ is linear.
On the other hand, when~$\Phi_0$ is not linear,
the free boundary condition in~\eqref{78:0ass12d}
takes into account the global behavior of the free boundary
and the (possibly fractional) perimeter of the minimizer
in the whole of the domain. In this sense,
this type of condition is ``self-driven'', since it is influenced
by the minimizer itself and not only by the environmental conditions
and the structural constants.

\section{Regularity of stationary points of the Alt-Caffarelli functional}
In this section we would like to discuss some recent results on the further connections 
between the Alt-Caffarelli problem and the minimal surfaces. 
More specifically, we consider the stationary points (in particular minimizers) of the functional 
\begin{equation}\label{78:APAJx:78}
{\mathscrboondox{E}}_{AC}[u]=\int_\Omega|\nabla u|^2+\lambda^2\chi_{\{u>0\}}\end{equation}
and the capillarity surfaces in the sphere of radius $\lambda$
(notice that the critical points of the functional in~\eqref{78:APAJx:78}
are related to the system
in~\eqref{78:APAJ}, see Theorem~2.5 in~\cite{MR618549} for details on this).

The starting point of our analysis is to study the classical capillary drop problem. 
\subsection{Capillary drop problem}
We first revisit the sessile drop problem and its higher dimensional analogue. 

\begin{figure}[h]  
\begin{center}
    \includegraphics[scale=0.85]{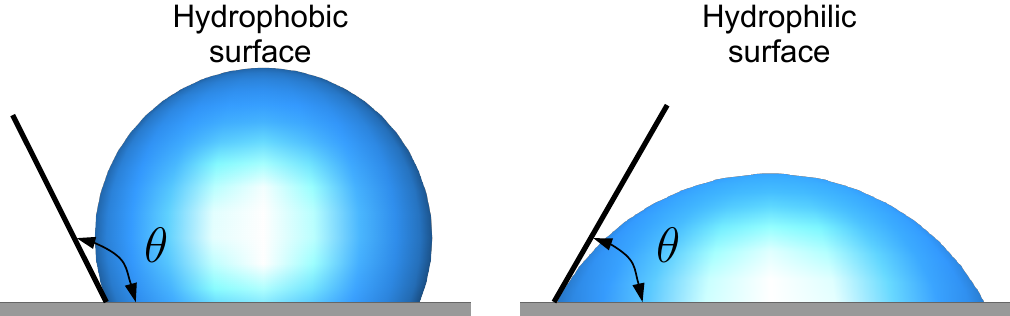}
  \end{center}
  \caption{Two regimes of drop-surface interactions.}\label{78534674325}
\end{figure}
We consider the functional 
\begin{equation}\label{sessile}
{\mathscrboondox{J}}(E):=\int_{\Omega_0}|D\chi_E|+g\int_{\Omega_0}x_{n+1}\chi_Edx-\int_{\partial \Omega_0}\lambda(x')\chi_Edx'\end{equation}
where 
$\Omega_0:=\{x=(x', x_{n+1}), x_{n+1}>0\}$, $x'=(x_1, \dots, x_n)\in \R^{n}$, 
$g>0$ is a given constant, $|\lambda(x')|<1$, and
$\chi_E$ is the characteristic function of $E\in \mathcal A$, where
$$\mathcal A:=\{E\subset \overline{ \Omega_0} \text{ s.t. } E \text{ has finite perimeter and } 
\mathcal H^{n+1}(E)=V\}.$$ Here, the parameter~$V>0$ is the volume fraction
of the droplet.

For $n=2$, the functional in~\eqref{sessile} is related to
the sessile drop problem, i.e. the problem of a 
(three-dimensional)
capillarity drop occupying the set $E$ and sitting in the
halfspace $\{x_{n+1}>0\}.$

We observe that the first term in ${\mathscrboondox{J}}(E)$ is the
energy due to the surface tension, 
the second term is the gravitational energy and the 
last term is the wetting energy which produces a contact angle 
$\theta(x')$ such that $\cos \theta(x')=\lambda(x')$
(see Figure~\ref{78534674325}).

By a Taylor expansion, we see that  
\[\sqrt{1+|\nabla u|^2}=1+\frac12|\nabla u|^2+\dots\]
Hence, if~$\partial E$ is the graph of a (smooth) function~$u\ge0$ (with small gradient),
we obtain the approximation
\begin{eqnarray*}&&\int_{\Omega_0}|D\chi_E|=
{\mathcal{H}}^{n-1}( (\partial E)\cap\Omega_0)=
\int_{(\partial\Omega_0)\cap\{u>0\}}\sqrt{1+|\nabla u|^2}\\&&\qquad=
\int_{\partial\Omega_0}\chi_{\{u>0\}}+\frac12
\int_{\partial\Omega_0}|\nabla u|^2+\dots\end{eqnarray*}
In other words, the functional ${\mathscrboondox{E}}_{AC}$ 
in~\eqref{78:APAJx:78}
is the linearization of the sessile drop problem
described by the functional~$
{\mathscrboondox{J}}$ in~\eqref{sessile}, with no gravity and constant
wetting energy density. This suggests
that there must be a strong link between the regularity of the minimizers of $
{\mathscrboondox{E}}_{AC}$ and the minimal surfaces.
We will now discuss\footnote{For completeness,
we recall that a nonlocal capillarity theory has been
recently developed in~\cite{2016arXiv160608610M, 2016arXiv161000075D}.} in which sense this link rigorously occurs.

\subsection{Homogeneity of blow-ups and the support function}
The first of such direct links was established in~\cite{MR2572253}, where it
is showed that the singular axisymmetric critical point of ${\mathscrboondox{E}}_{AC}$ is an energy minimizer in dimension 7. This
singular energy minimizer
of the Alt-Caffarelli problem can be seen as the analog of the Simons cone
\[S=\left\{x\in \R^8\ :\ \sum_{i=1}^4x_i^2=\sum_{i=5}^8x_i^2 \right\},\]
which is an example of a singular hypersurface of least perimeter in dimension~$ 8$. The minimality of the Simons cone was 
first proved by E. Bombieri, E. De Giorgi and E. Giusti in~\cite{MR0250205}.

The cones with non-negative mean curvature arise naturally in the blow-up procedure 
of the minimizer $u$ at a free boundary point. By Weiss' monotonicity formula (see~\cite{weiss}), 
any blow-up
limit $u_0$ of an energy minimizer of~${\mathscrboondox{E}}_{AC}$
is defined on $\R^n$ and
must be a homogeneous function of degree one.   

Let us write 
\begin{equation}\label{0oeyd8eu9refhido}
u_0(x)=rg(\sigma)\end{equation}
where $\sigma\in \mathbb S^{n-1}$ (being~$\mathbb S^{n-1}$, as usual,
the unit sphere in $\R^n$). Since $u_0$ is also a global minimizer of~$
{\mathscrboondox{E}}_{AC}$ then it follows that $\Delta u_0=0$ in $\Omega^+=\{u_0>0\}\cap \R^n$.
Rewriting the equation $\Delta u_0=0$ in polar coordinates we infer that~$g$
is a solution of the equation 
\begin{equation}\label{P098}
\Delta_{\mathbb S^{n-1}}g+n g=0.\end{equation}
 Here $\Delta_{\mathbb S^{n-1}}g$ is the Laplace-Beltrami operator
 on the sphere. 
 We observe that $u_0>0$ if and only if $g>0$ and $g=0$ on the free boundary of $u_0$ which is a 
 cone due to the homogeneity of $u_0$. 
   
Equation~\eqref{P098} can be rewritten as 
\begin{equation}\label{TRA}
\hbox{Trace}[g_{ij}+\delta_{ij}g]=0.\end{equation}
It is well-known that 
\begin{equation}\label{well}\begin{split}&
{\mbox{the eigenvalues of the matrix  $g_{ij}+\delta_{ij}g$}}\\&{\mbox{are the principal 
radii of curvature of the surface $S$}}\end{split}\end{equation} 
determined by the parameterization 
\[ {\mathbb S^{n-1}}\ni\sigma\longmapsto X(\sigma):=\sigma g(\sigma)+\nabla_\sigma g(\sigma),\]
see \cite{Aleksandrov}.
In addition, we have that $S$ and the sphere $\lambda \mathbb S^{n-1}$ are perpendicular at the 
contact points, see~\cite{AK}.

In this sense, one can interpret $g$ as the Minkowski support function of the surface $S.$ In other words
$X(\sigma)\cdot \sigma =g(\sigma)$ and it is the distance of the tangent plane with normal $\sigma$ from the origin.

\subsection{The mean radius equation}

The previous discussion tells us that the sum of the principal radii of the surface $S$
is zero. Indeed, let ${\kappa_i}=\frac1{R_i}, i=1, 2, \dots, n-1$ be the $i$th principal curvature 
of $S$ and $R_i$ the corresponding principal radius. Then, in view of~\eqref{well},
it holds that
the matrix~$g_{ij}+\delta_{ij}g$
has eigenvalues~$\frac1{\kappa_1},\dots,\frac1{\kappa_{n-1}}$,
and so its trace is equal to~$\sum_{i=1}^{n-1}\frac1{ \kappa_i}$.
{F}rom this and~\eqref{TRA}, we thus obtain that
\begin{equation}\label{RI}
\sum_{i=1}^{n-1}R_i=\sum_{i=1}^{n-1}\frac1{\kappa_i}=0\quad\hbox{in}\ \{g>0\}.\end{equation}
This is called the mean radius equation. 
Recalling~\eqref{0oeyd8eu9refhido},
the free boundary condition given in~\cite{MR618549} (and corresponding
to the constancy of~$|\nabla u_0|$) along~$\{g=0\}$
now becomes \[ |\nabla_\sigma g|=\lambda.\] This means that 
the surface $S$ is contained in the sphere of radius $\lambda$. 

We point out that in dimension~$n=2$ formula~\eqref{RI} reduces to 
\[\frac{\kappa_1+\kappa_2}{\kappa_1\kappa_2}=0\]
and therefore the mean curvature vanishes whenever the Gauss 
curvature is nonzero (i.e., $\kappa_1\kappa_2\not=0$). 
If $n\ge 3$, then such simple interpretation is not possible.

In terms of the classification of global cones for the Alt-Caffarelli problem, we recall that
the
open question is whether for $n=5, 6$ the stable solution of the mean radius equation 
in $\lambda \mathbb S^{n-1}$ such that $|\nabla_\sigma g|=\lambda$ on $g=0$ is
the disk passing through the origin (when~$n=2$, this question is settled in~\cite{MR618549},
the case~$n=3$ was addressed in~\cite{CJK} and~$n=4$ was proved in~\cite{JS}).

An open and challenging question is to classify the stationary solutions 
$g$ and 
the corresponding zero mean radius surfaces of given topological type.

\subsection{Flame models}
A closely related problem is the behavior of the solutions to 
the singular perturbation problem
\begin{equation}\label{pde-0}    
\left\{
\begin{array}{lll}
\Delta u_\e(x)=\beta_\e(u_\e) & \mbox{in}\quad B_1, \\
|u_\e|\le 1 & \mbox{in}\quad B_1,
\end{array}
\right.
\end{equation}
where $\e>0$ is small and
\begin{eqnarray}\label{beta}
\left\{\begin{array}{lll}
\beta_\e(t)=\frac1\e\beta\left(\frac{t}\e\right), \\
\,\\
\beta(t)\ge 0, \quad \supp\beta\subset [0, 1], \\
\,\\ \displaystyle\int_0^1\beta(t)dt= M>0
\end{array}
\right.
\end{eqnarray}
is an approximation of the Dirac measure.
It is well known that \eqref{pde-0}
models propagation
of equidiffusional premixed  flames with high activation of energy, see~\cite{MR1347971}.

The limit $u_0:=\lim\limits_{\e_j\to 0}u_{\e_j}$ (for a suitable sequence $\e_j\to 0$) solves a Bernoulli type free boundary problem
with the following free boundary condition
$$|\na u^+|^2-|\na u^-|^2=2M.$$
In fact, it holds that~$u_0$ is a stationary point of the Alt-Caffarelli functional
in~\eqref{78:APAJx:78}. 

If we choose $\{u_\e\}$ to be a family of minimizers of the functional
\begin{equation}
{\mathscrboondox{E}}_\e[u_\e]:=\int_{\Omega}\frac{|\na u_\e|^2}2
+\B(u_\e/\e), \quad{\mbox{ where }} \B(t)=\int_0^t\beta(s)ds,
\end{equation}
then $u_\e$ inherits the generic features of the
Alt-Caffarelli minimizers as described in~\cite{MR618549,
MR732100}
(e.g. non-degeneracy, rectifiability of $\fb u$, etc.).
Consequently, by
sending $\e\to 0$, one can see that the limit $u$ exists and it is a minimizer of
the Alt-Caffarelli functional $$ \int_{B_1}|\na u|^2+2M\chi_{\{u>0\}}.$$
As it was mentioned above, the singular set of minimizers is empty in dimensions~$
2, 3$ and  $4$. However, if $u_\e$ is not a minimizer then not much is known about the classification of the blow-ups of
the limit function $u$. An interesting question is to classify these stationary points of given topological type or Morse index. One recent result in this direction is that 
the associated surface $S$ that we constructed via the support function is of ring type then 
$S$ is a piece of catenoid, see \cite{AK}.

\def\cprime{$'$}

\vfill


\begin{thebibliography}{ADPM11}
\ifx \showCODEN  \undefined \def \showCODEN #1{CODEN #1}  \fi
\ifx \showISBN   \undefined \def \showISBN  #1{ISBN #1}   \fi
\ifx \showISSN   \undefined \def \showISSN  #1{ISSN #1}   \fi
\ifx \showLCCN   \undefined \def \showLCCN  #1{LCCN #1}   \fi
\ifx \showPRICE  \undefined \def \showPRICE #1{#1}        \fi
\ifx \showURL    \undefined \def \showURL {URL }          \fi
\ifx \path       \undefined \input path.sty               \fi
\ifx \ifshowURL \undefined
     \newif \ifshowURL
     \showURLtrue
\fi

\bibitem[AC81]{MR618549}
H.~W. Alt and L.~A. Caffarelli.
\newblock Existence and regularity for a minimum problem with free boundary.
\newblock {\em J. Reine Angew. Math.}, 325:\penalty0 105--144, 1981.
\newblock \showCODEN{JRMAA8}.
\newblock \showISSN{0075-4102}.

\bibitem[ACF84]{MR732100}
Hans~Wilhelm Alt, Luis~A. Caffarelli, and Avner Friedman.
\newblock Variational problems with two phases and their free boundaries.
\newblock {\em Trans. Amer. Math. Soc.}, 282\penalty0 (2):\penalty0 431--461,
  1984.
\newblock \showCODEN{TAMTAM}.
\newblock \showISSN{0002-9947}.
\newblock \ifshowURL {\showURL \path|http://dx.doi.org/10.2307/1999245|}\fi.

\bibitem[ACKS01]{MR1808651}
I.~Athanasopoulos, L.~A. Caffarelli, C.~Kenig, and S.~Salsa.
\newblock An area-{D}irichlet integral minimization problem.
\newblock {\em Comm. Pure Appl. Math.}, 54\penalty0 (4):\penalty0 479--499,
  2001.
\newblock \showCODEN{CPAMA}.
\newblock \showISSN{0010-3640}.
\newblock \ifshowURL {\showURL
  \path|http://dx.doi.org/10.1002/1097-0312(200104)54:4<479::AID-CPA3>3.3.CO;2-U|}\fi.

\bibitem[ADPM11]{MR2765717}
Luigi Ambrosio, Guido De~Philippis, and Luca Martinazzi.
\newblock Gamma-convergence of nonlocal perimeter functionals.
\newblock {\em Manuscripta Math.}, 134\penalty0 (3-4):\penalty0 377--403, 2011.
\newblock \showCODEN{MSMHB2}.
\newblock \showISSN{0025-2611}.
\newblock \ifshowURL {\showURL
  \path|http://dx.doi.org/10.1007/s00229-010-0399-4|}\fi.

\bibitem[Ale39]{Aleksandrov}
A.~Alexandroff.
\newblock \"uber die {O}berfl\"achenfunktion eines konvexen {K}\"orpers.
  ({B}emerkung zur {A}rbeit ``{Z}ur {T}heorie der gemischten {V}olumina von
  konvexen {K}\"orpern'').
\newblock {\em Rec. Math. N.S. [Mat. Sbornik]}, 6(48):\penalty0 167--174, 1939.

\bibitem[AV14]{MR3230079}
Nicola Abatangelo and Enrico Valdinoci.
\newblock A notion of nonlocal curvature.
\newblock {\em Numer. Funct. Anal. Optim.}, 35\penalty0 (7-9):\penalty0
  793--815, 2014.
\newblock \showISSN{0163-0563}.
\newblock \ifshowURL {\showURL
  \path|http://dx.doi.org/10.1080/01630563.2014.901837|}\fi.

\bibitem[BBM02]{MR1945278}
Jean Bourgain, Ha{\"{\i}}m Brezis, and Petru Mironescu.
\newblock Limiting embedding theorems for {$W^{s,p}$} when {$s\uparrow1$} and
  applications.
\newblock {\em J. Anal. Math.}, 87:\penalty0 77--101, 2002.
\newblock \showCODEN{JOAMAV}.
\newblock \showISSN{0021-7670}.
\newblock \ifshowURL {\showURL \path|http://dx.doi.org/10.1007/BF02868470|}\fi.
\newblock Dedicated to the memory of Thomas H. Wolff.

\bibitem[BDGG69]{MR0250205}
E.~Bombieri, E.~De~Giorgi, and E.~Giusti.
\newblock Minimal cones and the {B}ernstein problem.
\newblock {\em Invent. Math.}, 7:\penalty0 243--268, 1969.
\newblock \showISSN{0020-9910}.
\newblock \ifshowURL {\showURL \path|http://dx.doi.org/10.1007/BF01404309|}\fi.

\bibitem[BFV14]{MR3331523}
Bego{\~n}a Barrios, Alessio Figalli, and Enrico Valdinoci.
\newblock Bootstrap regularity for integro-differential operators and its
  application to nonlocal minimal surfaces.
\newblock {\em Ann. Sc. Norm. Super. Pisa Cl. Sci. (5)}, 13\penalty0
  (3):\penalty0 609--639, 2014.
\newblock \showISSN{0391-173X}.

\bibitem[BLV16]{2016arXiv161208295B}
Claudia {Bucur}, Luca {Lombardini}, and Enrico {Valdinoci}.
\newblock {Complete stickiness of nonlocal minimal surfaces for small values of
  the fractional parameter}.
\newblock {\em ArXiv e-prints}, December 2016.

\bibitem[BN16]{MR3556344}
Ha{\"{\i}}m Brezis and Hoai-Minh Nguyen.
\newblock The {BBM} formula revisited.
\newblock {\em Atti Accad. Naz. Lincei Rend. Lincei Mat. Appl.}, 27\penalty0
  (4):\penalty0 515--533, 2016.
\newblock \showISSN{1120-6330}.
\newblock \ifshowURL {\showURL \path|http://dx.doi.org/10.4171/RLM/746|}\fi.

\bibitem[Caf95]{MR1347971}
Luis~A. Caffarelli.
\newblock Uniform {L}ipschitz regularity of a singular perturbation problem.
\newblock {\em Differential Integral Equations}, 8\penalty0 (7):\penalty0
  1585--1590, 1995.
\newblock \showISSN{0893-4983}.

\bibitem[CDSS16]{MR3532394}
L.~Caffarelli, D.~De~Silva, and O.~Savin.
\newblock Obstacle-type problems for minimal surfaces.
\newblock {\em Comm. Partial Differential Equations}, 41\penalty0 (8):\penalty0
  1303--1323, 2016.
\newblock \showCODEN{CPDIDZ}.
\newblock \showISSN{0360-5302}.
\newblock \ifshowURL {\showURL
  \path|http://dx.doi.org/10.1080/03605302.2016.1192646|}\fi.

\bibitem[CJK04]{CJK}
Luis~A. Caffarelli, David Jerison, and Carlos~E. Kenig.
\newblock Global energy minimizers for free boundary problems and full
  regularity in three dimensions.
\newblock In {\em Noncompact problems at the intersection of geometry,
  analysis, and topology}, volume 350 of {\em Contemp. Math.}, pages 83--97.
  Amer. Math. Soc., Providence, RI, 2004.
\newblock \ifshowURL {\showURL
  \path|http://dx.doi.org/10.1090/conm/350/06339|}\fi.

\bibitem[CRS10a]{MR2675483}
L.~Caffarelli, J.-M. Roquejoffre, and O.~Savin.
\newblock Nonlocal minimal surfaces.
\newblock {\em Comm. Pure Appl. Math.}, 63\penalty0 (9):\penalty0 1111--1144,
  2010.
\newblock \showCODEN{CPAMA}.
\newblock \showISSN{0010-3640}.
\newblock \ifshowURL {\showURL \path|http://dx.doi.org/10.1002/cpa.20331|}\fi.

\bibitem[CRS10b]{MR2677613}
Luis~A. Caffarelli, Jean-Michel Roquejoffre, and Yannick Sire.
\newblock Variational problems for free boundaries for the fractional
  {L}aplacian.
\newblock {\em J. Eur. Math. Soc. (JEMS)}, 12\penalty0 (5):\penalty0
  1151--1179, 2010.
\newblock \showISSN{1435-9855}.
\newblock \ifshowURL {\showURL \path|http://dx.doi.org/10.4171/JEMS/226|}\fi.

\bibitem[CSV15]{MR3390089}
Luis Caffarelli, Ovidiu Savin, and Enrico Valdinoci.
\newblock Minimization of a fractional perimeter-{D}irichlet integral
  functional.
\newblock {\em Ann. Inst. H. Poincar\'e Anal. Non Lin\'eaire}, 32\penalty0
  (4):\penalty0 901--924, 2015.
\newblock \showISSN{0294-1449}.
\newblock \ifshowURL {\showURL
  \path|http://dx.doi.org/10.1016/j.anihpc.2014.04.004|}\fi.

\bibitem[CV11]{MR2782803}
Luis Caffarelli and Enrico Valdinoci.
\newblock Uniform estimates and limiting arguments for nonlocal minimal
  surfaces.
\newblock {\em Calc. Var. Partial Differential Equations}, 41\penalty0
  (1-2):\penalty0 203--240, 2011.
\newblock \showISSN{0944-2669}.
\newblock \ifshowURL {\showURL
  \path|http://dx.doi.org/10.1007/s00526-010-0359-6|}\fi.

\bibitem[CV13]{MR3107529}
Luis Caffarelli and Enrico Valdinoci.
\newblock Regularity properties of nonlocal minimal surfaces via limiting
  arguments.
\newblock {\em Adv. Math.}, 248:\penalty0 843--871, 2013.
\newblock \showISSN{0001-8708}.
\newblock \ifshowURL {\showURL
  \path|http://dx.doi.org/10.1016/j.aim.2013.08.007|}\fi.

\bibitem[D{\'a}v02]{MR1942130}
J.~D{\'a}vila.
\newblock On an open question about functions of bounded variation.
\newblock {\em Calc. Var. Partial Differential Equations}, 15\penalty0
  (4):\penalty0 519--527, 2002.
\newblock \showISSN{0944-2669}.
\newblock \ifshowURL {\showURL
  \path|http://dx.doi.org/10.1007/s005260100135|}\fi.

\bibitem[DFPV13]{MR3007726}
Serena Dipierro, Alessio Figalli, Giampiero Palatucci, and Enrico Valdinoci.
\newblock Asymptotics of the {$s$}-perimeter as {$s\searrow0$}.
\newblock {\em Discrete Contin. Dyn. Syst.}, 33\penalty0 (7):\penalty0
  2777--2790, 2013.
\newblock \showISSN{1078-0947}.
\newblock \ifshowURL {\showURL
  \path|http://dx.doi.org/10.3934/dcds.2013.33.2777|}\fi.

\bibitem[Dip14]{MR3506705}
S.~Dipierro.
\newblock Asymptotics of fractional perimeter functionals and related problems.
\newblock {\em Rend. Semin. Mat. Univ. Politec. Torino}, 72\penalty0
  (1-2):\penalty0 3--16, 2014.
\newblock \showISSN{0373-1243}.

\bibitem[DK15]{2015arXiv150807447D}
Serena {Dipierro} and Aram~L. {Karakhanyan}.
\newblock {Stratification of free boundary points for a two-phase variational
  problem}.
\newblock {\em ArXiv e-prints}, August 2015.

\bibitem[DKV15]{2015arXiv151203043D}
Serena {Dipierro}, Aram {Karakhanyan}, and Enrico {Valdinoci}.
\newblock {A class of unstable free boundary problems}.
\newblock {\em ArXiv e-prints}, December 2015.

\bibitem[DKV16]{2016arXiv161100412D}
Serena {Dipierro}, Aram {Karakhanyan}, and Enrico {Valdinoci}.
\newblock {A nonlinear free boundary problem with a self-driven Bernoulli
  condition}.
\newblock {\em ArXiv e-prints}, November 2016.

\bibitem[DLV]{LUCA}
Serena Dipierro, Luca Lombardini, and Enrico Valdinoci.
\newblock A free boundary problem: superposition of nonlocal energy plus
  classical perimeter.
\newblock {\em {P}reprint}.

\bibitem[DMV16]{2016arXiv161000075D}
Serena {Dipierro}, Francesco {Maggi}, and Enrico {Valdinoci}.
\newblock {Asymptotic expansions of the contact angle in nonlocal capillarity
  problems}.
\newblock {\em ArXiv e-prints}, September 2016.

\bibitem[DNPV12]{MR2944369}
Eleonora Di~Nezza, Giampiero Palatucci, and Enrico Valdinoci.
\newblock Hitchhiker's guide to the fractional {S}obolev spaces.
\newblock {\em Bull. Sci. Math.}, 136\penalty0 (5):\penalty0 521--573, 2012.
\newblock \showISSN{0007-4497}.
\newblock \ifshowURL {\showURL
  \path|http://dx.doi.org/10.1016/j.bulsci.2011.12.004|}\fi.

\bibitem[DSJ09]{MR2572253}
Daniela De~Silva and David Jerison.
\newblock A singular energy minimizing free boundary.
\newblock {\em J. Reine Angew. Math.}, 635:\penalty0 1--21, 2009.
\newblock \showISSN{0075-4102}.
\newblock \ifshowURL {\showURL
  \path|http://dx.doi.org/10.1515/CRELLE.2009.074|}\fi.

\bibitem[DSR12]{MR2926238}
D.~De~Silva and J.~M. Roquejoffre.
\newblock Regularity in a one-phase free boundary problem for the fractional
  {L}aplacian.
\newblock {\em Ann. Inst. H. Poincar\'e Anal. Non Lin\'eaire}, 29\penalty0
  (3):\penalty0 335--367, 2012.
\newblock \showISSN{0294-1449}.
\newblock \ifshowURL {\showURL
  \path|http://dx.doi.org/10.1016/j.anihpc.2011.11.003|}\fi.

\bibitem[DSV15]{MR3427047}
Serena Dipierro, Ovidiu Savin, and Enrico Valdinoci.
\newblock A nonlocal free boundary problem.
\newblock {\em SIAM J. Math. Anal.}, 47\penalty0 (6):\penalty0 4559--4605,
  2015.
\newblock \showISSN{0036-1410}.
\newblock \ifshowURL {\showURL \path|http://dx.doi.org/10.1137/140999712|}\fi.

\bibitem[DSV16a]{Dipierro2016JFA}
Serena Dipierro, Ovidiu Savin, and Enrico Valdinoci.
\newblock Boundary behavior of nonlocal minimal surfaces.
\newblock {\em J. Funct. Anal.}, 2016.
\newblock \showISSN{0022-1236}.
\newblock \ifshowURL {\showURL
  \path|http://www.sciencedirect.com/science/article/pii/S0022123616303858|}\fi.

\bibitem[DSV16b]{MR3516886}
Serena Dipierro, Ovidiu Savin, and Enrico Valdinoci.
\newblock Graph properties for nonlocal minimal surfaces.
\newblock {\em Calc. Var. Partial Differential Equations}, 55\penalty0
  (4):\penalty0 Paper No. 86, 25, 2016.
\newblock \showISSN{0944-2669}.
\newblock \ifshowURL {\showURL
  \path|http://dx.doi.org/10.1007/s00526-016-1020-9|}\fi.

\bibitem[DV15]{MR3320130}
Serena Dipierro and Enrico Valdinoci.
\newblock On a fractional harmonic replacement.
\newblock {\em Discrete Contin. Dyn. Syst.}, 35\penalty0 (8):\penalty0
  3377--3392, 2015.
\newblock \showISSN{1078-0947}.
\newblock \ifshowURL {\showURL
  \path|http://dx.doi.org/10.3934/dcds.2015.35.3377|}\fi.

\bibitem[DV16]{Dipierro2016HP}
Serena Dipierro and Enrico Valdinoci.
\newblock Continuity and density results for a one-phase nonlocal free boundary
  problem.
\newblock {\em Ann. Inst. H. Poincar\'e Anal. Non Lin\'eaire}, 2016.
\newblock \showISSN{0294-1449}.
\newblock \ifshowURL {\showURL
  \path|http://www.sciencedirect.com/science/article/pii/S0294144916300853|}\fi.

\bibitem[DV17]{2016arXiv160706872D}
Serena {Dipierro} and Enrico {Valdinoci}.
\newblock Nonlocal minimal surfaces: interior regularity, quantitative
  estimates and boundary stickiness.
\newblock {\em Recent Developments in the Nonlocal Theory (T. Kuusi, G.
  Palatucci eds.). Book Series on Measure Theory. De Gruyter, Berlin}, 2017.

\bibitem[FV15]{FIGALLI}
Alessio Figalli and Enrico Valdinoci.
\newblock Regularity and bernstein-type results for nonlocal minimal surfaces.
\newblock {\em J. Reine Angew. Math.}, 2015.
\newblock \showISSN{0075-4102}.
\newblock \ifshowURL {\showURL
  \path|https://www.degruyter.com/view/j/crll.ahead-of-print/crelle-2015-0006/crelle-2015-0006.xml|}\fi.

\bibitem[JS15]{JS}
David Jerison and Ovidiu Savin.
\newblock Some remarks on stability of cones for the one-phase free boundary
  problem.
\newblock {\em Geom. Funct. Anal.}, 25\penalty0 (4):\penalty0 1240--1257, 2015.
\newblock \showISSN{1016-443X}.
\newblock \ifshowURL {\showURL
  \path|http://dx.doi.org/10.1007/s00039-015-0335-6|}\fi.

\bibitem[Kar16]{AK}
Aram Karakhanyan.
\newblock Minimal surfaces arising in singular perturbation problems.
\newblock {\em {P}reprint}, 2016.

\bibitem[MS02]{MR1940355}
V.~Maz{\cprime}ya and T.~Shaposhnikova.
\newblock On the {B}ourgain, {B}rezis, and {M}ironescu theorem concerning
  limiting embeddings of fractional {S}obolev spaces.
\newblock {\em J. Funct. Anal.}, 195\penalty0 (2):\penalty0 230--238, 2002.
\newblock \showCODEN{JFUAAW}.
\newblock \showISSN{0022-1236}.
\newblock \ifshowURL {\showURL
  \path|http://dx.doi.org/10.1006/jfan.2002.3955|}\fi.

\bibitem[MV16]{2016arXiv160608610M}
Francesco {Maggi} and Enrico {Valdinoci}.
\newblock {Capillarity problems with nonlocal surface tension energies}.
\newblock {\em ArXiv e-prints}, June 2016.

\bibitem[Pon04]{MR2033060}
Augusto~C. Ponce.
\newblock A new approach to {S}obolev spaces and connections to
  {$\Gamma$}-convergence.
\newblock {\em Calc. Var. Partial Differential Equations}, 19\penalty0
  (3):\penalty0 229--255, 2004.
\newblock \showISSN{0944-2669}.
\newblock \ifshowURL {\showURL
  \path|http://dx.doi.org/10.1007/s00526-003-0195-z|}\fi.

\bibitem[Ser71]{MR0333220}
James Serrin.
\newblock A symmetry problem in potential theory.
\newblock {\em Arch. Rational Mech. Anal.}, 43:\penalty0 304--318, 1971.
\newblock \showISSN{0003-9527}.

\bibitem[SV13a]{MR3090533}
Ovidiu Savin and Enrico Valdinoci.
\newblock Regularity of nonlocal minimal cones in dimension 2.
\newblock {\em Calc. Var. Partial Differential Equations}, 48\penalty0
  (1-2):\penalty0 33--39, 2013.
\newblock \showISSN{0944-2669}.
\newblock \ifshowURL {\showURL
  \path|http://dx.doi.org/10.1007/s00526-012-0539-7|}\fi.

\bibitem[SV13b]{MR3035063}
Ovidiu Savin and Enrico Valdinoci.
\newblock Some monotonicity results for minimizers in the calculus of
  variations.
\newblock {\em J. Funct. Anal.}, 264\penalty0 (10):\penalty0 2469--2496, 2013.
\newblock \showCODEN{JFUAAW}.
\newblock \showISSN{0022-1236}.
\newblock \ifshowURL {\showURL
  \path|http://dx.doi.org/10.1016/j.jfa.2013.02.005|}\fi.

\bibitem[Wei98]{weiss}
Georg~S. Weiss.
\newblock Partial regularity for weak solutions of an elliptic free boundary
  problem.
\newblock {\em Comm. Partial Differential Equations}, 23\penalty0
  (3-4):\penalty0 439--455, 1998.
\newblock \showISSN{0360-5302}.
\newblock \ifshowURL {\showURL
  \path|http://dx.doi.org/10.1080/03605309808821352|}\fi.

\end{thebibliography}
\end{document}